\documentclass{amsart}
\numberwithin{equation}{section}
\usepackage{graphicx}
\usepackage{wrapfig}
\newtheorem{theorem}{Theorem}[section]
\newtheorem{lemma}[theorem]{Lemma}
\newtheorem{corollary}[theorem]{Corollary}

\newtheorem{definition}[theorem]{Definition}

\begin{document}

\author{Ravshan Ashurov}
\address{Ashurov R: V. I. Romanovskiy Institute of Mathematics, Uzbekistan Academy of Sciences, Tashkent, Uzbekistan}
\email{ashurovr@gmail.com}

\author{Ilyoskhuja Sulaymonov}
\address{Sulaymonov I: National University of Uzbekistan,  Tashkent, Uzbekistan}
\email{ilyosxojasulaymonov@gmail.com}

\small

\title[Monotonicity in the parameter of the Mittag-Leffler function ...] {Monotonicity in the parameter of the Mittag-Leffler function and determining the fractional exponent of the subdiffusion equation}

\begin{abstract}
In this paper, we prove the strict monotonicity in the parameter $\rho$ of the Mittag-Leffler functions $E_\rho(-t^\rho)$ and $t^{\rho -1}E_{\rho,\rho}(-t^\rho)$. Then, these results are applied to solve the inverse problem of determining the order of the fractional derivative in subdiffusion equations, where the available measurement is given at one point in space-time. In particular, we find the missing conditions in the previously known work in this area. Moreover, the obtained results are valid for a wider class of subdiffusion equations than those considered previously. An example of an initial boundary value problem constructed by Sh.A. Alimov is given, for which the inverse problem under consideration has a unique solution. We also point out the application of the monotonicity of the Mittag-Leffler functions to solving some other inverse problems of determining the order of a fractional derivative.

\end{abstract}

\maketitle

\section{Introduction}

Let $\Omega$ be an arbitrary $N$ - dimensional domain with a sufficiently smooth boundary $\partial\Omega$. Consider the following \textit{Initial-boundary value problem}:
\begin{equation}\label{prob1}
\left\{
\begin{aligned}
&D_t^\rho u(x,t) -\Delta u(x,t) =0,\quad x\in\Omega,\quad 0< t\leq T,\\
&u{(x,t)}|_{\partial\Omega}=0, \\
&u(x,0)=\varphi(x),\quad x\in\Omega,
\end{aligned}
\right.
\end{equation}
 where $\varphi(x)$ is a continuous function, $\Delta$ is the Laplace operator, $\rho\in(0,1)$ and $D_t^\rho$ is the fractional Caputo derivative defined as (see, for example, \cite{Klb}, p. 91):
$$
D_t^\rho u(x,t)=\frac{1}{\Gamma
(1-\rho)}\frac{d}{dt}\int\limits_0^t\frac{u(x,\xi)-u(x,0)}{(t-\xi)^{\rho}} d\xi,
\quad t>0,\,\, x\in\Omega.
$$

When modeling various processes, the order of the fractional derivative $\rho$ is often unknown. Unfortunately, in the processes under consideration there is no device for measuring this parameter. One of the effective methods for finding the unknown order of the derivative is, by setting some additional conditions for a solution $u(x,t)$ of problem (\ref{prob1}), to solve analytically the corresponding inverse problem.

In recent years, significant advances have been made in the study of this type of inverse problems (see Liu et al. \cite{Liu} and references therein for a comprehensive review of the work published up to 2019). Moreover, numerous studies (e.g. \cite{AU}-\cite{AM}) published after 2019 have proposed additional conditions for efficiently treating and solving these inverse problems.

A very interesting paper \cite{GongLi}  by G. Li, Z. Wang, X. Jia, Y. Zhang was recently published, in which the order of the fractional derivative $\rho$ in the equation $D_t^\rho u(x,t) - c \Delta u (x,t) = 0$ ($\Omega \subset \mathbb{R}^d, d=1,2,3$ and $c$ is a sufficiently small positive number) was reconstructed from information on the solution of the initial-boundary value problem at only one point of space-time $(x_0, t_0)$ (see Theorem 2). However, in our opinion, this result is questionable in the sense that the conditions of this theorem are not sufficient to guarantee its assertion. Note also that the proof relies on the condition $d\le 3$. However, the paper contains an interesting auxiliary result (see Theorem 1): the authors proved that the function $g(\rho) = E_\rho(-ct^\rho)$ has a negative derivative, $g'(\rho) < 0$, if $t$ is large enough and $c$ is small enough.

Motivated by the results and interesting ideas of \cite{GongLi}, we continued our research in this direction. We managed to find the missing condition in Theorem 2 and remove some restrictions mentioned in \cite{GongLi}. In particular:

1) we have proved that the Mittag-Leffler functions $t^{\rho -1}E_{\rho,\rho}(-t^\rho)$ and $E_\rho(-t^\rho)$ increase monotonically in parameter $\rho$ for sufficiently small $t$. This is the main result of this work;

2) the requirement of a sufficiently small coefficient $c$ was removed from the considered equation;

3) the restriction on the dimension of the domain $\Omega \subset \mathbb{R}^d, d=1,2,3$ was removed. Of course, the transition to the domain $\Omega \subset \mathbb{R}^N$, with $N\geq 1$, does not cause any particular difficulty, but we would like to draw the attention of readers to one fundamental work by V.A. Ilyin \cite{Il};

4) it is shown that the conditions of Theorem 2 of \cite{GongLi} is not complete and an updated version of this theorem with a proof is presented.

5) other applications of the monotonicity in parameter of Mittag-Leffler functions are given.

The remainder of the paper is structured as follows.

Section 2 introduces the necessary auxiliary concepts. Section 3 proves the monotonicity of the Mittag-Leffler functions in the parameter $\rho$. This result, which is naturally of independent interest, can be applied to the solution of various inverse problems of determining the order of a fractional derivative; some applications are noted at the end of this section. Section 4 studies the initial-boundary value problem (\ref{prob1}). Section 5 discusses the inverse problem of determining the fractional order using information on the solution of the initial-boundary value problem at only one point of space-time $(x_0, t_0)$. This requires rather strict conditions on the corresponding eigenfunctions and Fourier coefficients of the initial function. In the next section an example of an initial boundary value problem constructed by Sh.A. Alimov is given, for which this inverse problem has always a unique solution.  In the final section, using the monotonicity of the Mittag-Leffler function, a new proof of one result of A.V. Pskhu is presented.

\section{Preliminaries}

In this section we recall the fundamental result of V.A. Ilyin \cite{Il} on the convergence of Fourier coefficients and present an estimate of the Mittag-Leffler function.

Let $\Omega$ be a bounded $N$-dimentional domain with a sufficiently smooth boundary $\partial\Omega$ and $\{v_k(x)\}$ denote the complete system of orthonormal in $L_2(\Omega)$  eigenfunctions and $\{\lambda_k\}$ the set of positive eigenvalues of the spectral problem:
\[
\left\{
\begin{aligned}
& -\Delta v(x) = \lambda v(x), \quad x \in \Omega, \\
& v(x) \big|_{\partial \Omega} = 0.
\end{aligned}
\right.
\]
Let us present some assertions about eigenfunctions $v_k(x)$ and eigenvalues $\lambda_k$ proved by V.A.Ilyin \cite{Il}.
\begin{lemma}\label{Illem1} The series $\sum\limits_{k=1}^\infty\lambda_k^{-\left(\left[\frac{N}{2}\right]+1\right)}v_k^2(x)$ converges uniformly in a closed domain $\overline{\Omega}$.
\end{lemma}
\begin{lemma}\label{Illem2} Let the function $g(x)$ satisfy the conditions
\begin{enumerate}
    \item $g(x)\in C^p(\overline{\Omega}),\,\, \frac{\partial^{p+1}g(x)}{\partial x_1^{p_1}\dots \partial x_n^{p_n}}\in L_2(\Omega),\,\, p+1=p_1+p_2+\dots p_n, p\ge1,$
    \item $g(x)|_{\partial\Omega}=\Delta g(x)|_{\partial\Omega}=\dots=\Delta^{\left[\frac{p}{2}\right]}g(x)|_{\partial\Omega}=0.$
\end{enumerate}
    Then the number series $\sum\limits_{k=1}^\infty g_k^2\lambda_k^{p+1}$ converges, where $g_k=(g,v_k)$.
\end{lemma}

Let $A$ be a self-adjoint extension in $L_2(\Omega)$ of the Laplace operator with the homogeneous Dirichlet condition. Let $v_k(x)$ be the eigenfunctions and $\lambda_k$ the corresponding eigenvalues of this operator. For any real number $\tau$ we define the power of the operator $A$ by the formula
\[
A^\tau f(x) = \sum\limits_{k=1}^\infty \lambda_k^\tau f_k v_k(x), \quad f_k=(f, v_k),
\]
the domain of definition has the form
\[
D(A^\tau)= \{ f\in L_2(\Omega), \quad \sum\limits_{k=1}^\infty \lambda_k^{2\tau} |f_k|^2<\infty \}.
\]

\begin{definition}
Let $M$ be a fixed number. We will say that $f\in I_M(\Omega)$ if
\begin{enumerate}
    \item satisfy conditions Lemma \ref{Illem2} with the exponent $p=\left[\frac{N}{2}\right]$
    \item $||A^{\frac{1}{2}\left([\frac{N}{2}]+1\right)}f||_{L_2(\Omega)} \leq M.$
\end{enumerate}
\end{definition}

For $0 < \rho < 1$, let $E_{\rho,\mu}(z)$ denote the Mittag-Leffler function defined as:
$$
E_{\rho,\mu}(z)= \sum\limits_{k=0}^\infty \frac{z^k}{\Gamma(\rho
k+\mu)},\quad \mu, z\in \mathbb{C}.
$$
If \(\mu = 1\), then the Mittag-Leffler function is called the one-parameter or classical Mittag-Leffler function and is denoted by \(E_{\rho}(z) = E_{\rho,1}(z)\).

Recall the following estimate of the Mittag-Leffler functions (see, e.g. \cite{Gor}, p. 29).

\begin{lemma}\label{mll4} For any $t\geq 0$ one has
$$
|E_{\rho,\mu}(-t)|\leq \frac{C}{1+t},\quad \mu\in \mathbb{C},
$$
where constant $C$ does not depend on $t$.
\end{lemma}

\section{Monotonicity of the Mittag-Leffler functions}

This paragraph uses some original ideas from the paper \cite{GongLi}  by G. Li, Z. Wang, X. Jia, Y. Zhang.

Let us first introduce some important concepts to prove the strict monotonicity of the Mittag-Leffler function $E_{\rho}(-t^\rho)$.
\begin{lemma}\label{lem1}
Let $\rho\in(0,1)$. Then the following equality holds:
$$
\lim\limits_{n\rightarrow\infty}\frac{\Gamma(\rho n)}{\Gamma(\rho n+\rho)}=0.
$$
\end{lemma}
This lemma was proved in \cite{GongLi}. For the convenience of the reader, we recall this proof.

\begin{proof} Using the Stirling approximation of the gamma function we have:
$$
\frac{\Gamma(\rho n)}{\Gamma(\rho n+\rho)}\sim \frac{\frac{\sqrt{2\pi}}{\sqrt{\rho n}}\left(\frac{\rho n}{e}\right)^{\rho n}}{\frac{\sqrt{2\pi}}{\sqrt{\rho n+\rho}}\left(\frac{\rho n+\rho}{e}\right)^{\rho n+\rho}}\sim e^{\rho}\left(\frac{\rho n}{\rho n+\rho}\right)^{\rho n}\frac{1}{(\rho n+\rho)^\rho}
$$
$$
\times\sqrt{\frac{\rho n+\rho}{\rho n}}
\sim e^{\rho} \left[\left(\frac{1}{1+1/n}\right)^n\right]^\rho \frac{1}{(\rho n+\rho)^\rho}\sqrt{1+\frac{1}{n}} \quad (n\rightarrow \infty).
$$
According to $\lim\limits_{n\rightarrow\infty} \left(\frac{1}{1+1/n}\right)^n=e^{-1}$, we get:
$$
\lim\limits_{n\rightarrow\infty} \frac{\Gamma(\rho n)}{\Gamma(\rho n+\rho)}=\lim\limits_{n\rightarrow\infty} e^\rho e^{-\rho} \frac{1}{(\rho n+\rho)^\rho}\sqrt{1+\frac{1}{n}}=\lim\limits_{n\rightarrow\infty}  \frac{1}{(\rho n+\rho)^\rho}=0.
$$

Lemma \ref{lem1} is proved.
\end{proof}

Now we calculate the derivative of the Mittag-Leffler function $E_{\rho}(-t^\rho)$. Let $\Phi(\rho)$ be the logarithmic derivative of the gamma function $\Gamma(\rho)$ (see, for example \cite{Bateman}). Then $\Gamma'(\rho)=\Gamma(\rho)\Phi(\rho)$ and therefore, we have:
$$
\frac{d}{d\rho}E_{\rho}(-t^\rho)=\frac{d}{d\rho}\sum\limits_{n=0}^\infty (-1)^n\frac{t^{\rho n}}{\Gamma(\rho n+1)}
$$
\begin{equation}\label{derML}
=\sum\limits_{n=1}^\infty (-1)^nnt^{\rho n}\frac{\ln{t}-\Phi(\rho n+1)}{\Gamma(\rho n+1)}=\sum\limits_{n=1}^\infty  (-1)^ny_n,
\end{equation}
where
$$
\Phi(\rho n+1)=-\gamma-\frac{1}{\rho n+1}+\sum\limits_{s=1}^\infty\bigg(\frac{1}{s}-\frac{1}{s+\rho n+1}\bigg),
  $$
\begin{equation}\label{sery_n}
y_n=nt^{\rho n}\frac{\ln{t}-\Phi(\rho n+1)}{\Gamma(\rho n+1)},
\end{equation}
and $\gamma\approx 0.57722$ is the Euler-Mascheroni constant.

\begin{lemma}\label{lem2}
 Let $\rho\in(0,1)$. Then the following two statements hold:

1)
  $$
\lim\limits_{n\rightarrow\infty} \frac{\Phi(\rho n+\rho+1)}{\Phi(\rho n+1)}=1.
  $$

2) $\Phi(\rho n+1)> -\frac{1}{\rho +1}$ for all $n\ge1$.
\end{lemma}

\begin{proof} Part 1 of this lemma was proved in \cite{GongLi}. We give a slightly different proof.
Consider the function $f(z)=z(\ln{z}-\Phi(z))$. Function $f(z)$ is strictly decreasing and strictly convex on $(0, \infty)$ (see, \cite{Anderson}) and the following relationships are fulfilled:
\begin{equation}\label{f(z)lim}
    \lim\limits_{z\rightarrow0}f(z)=1 \quad \text{and} \quad \lim\limits_{z\rightarrow\infty}f(z)=\frac{1}{2}.
\end{equation}
From (\ref{f(z)lim}) and the monotonicity of $f(z)$ we get:
\begin{equation}\label{estfz}
\frac{1}{2}\le f(z)\le 1, \quad z\in(0,\infty).
\end{equation}
Using estimate (\ref{estfz}) and definition of function $f(z)$ we obtain
\begin{equation}\label{estpsi}
    \ln{z}-\frac{1}{z}\le \Phi(z)\le \ln{z}-\frac{1}{2z}.
\end{equation}
Apply estimate (\ref{estpsi}) to have:
\begin{equation}\label{estgam}
\ln{(\rho n+1)}-\frac{1}{\rho n+1}\le \Phi(\rho n+1)\le \ln{(\rho n+1)}-\frac{1}{2(\rho n+1)}.
\end{equation}
Therefore,
$$
\frac{\ln{(\rho n+\rho+1)}-\frac{1}{\rho n+\rho+1}}{\ln{(\rho n+1)}-\frac{1}{2(\rho n+1)}}\le\frac{\Phi(\rho n+\rho+1)}{\Phi(\rho n+1)}\le \frac{\ln{(\rho n+\rho+1)}-\frac{1}{2(\rho n+\rho+1)}}{\ln{(\rho n+1)}-\frac{1}{\rho n+1}}.
$$
Firstly we will show that
\begin{equation}\label{ln/ln}
\lim\limits_{n\rightarrow\infty} \frac{\ln{(\rho n+\rho+1)}}{\ln{(\rho n+1)}}=1.
\end{equation}
Using Lopital's theorem we get
$$
\lim\limits_{n\rightarrow\infty} \frac{\ln{(\rho n+\rho+1)}}{\ln{(\rho n+1)}}=\lim\limits_{n\rightarrow\infty} \frac{\frac{\rho}{\rho n+\rho+1}}{\frac{\rho}{\rho n+1}}=\lim\limits_{n\rightarrow\infty} \frac{\rho n+1}{\rho n+\rho+1}=1
$$
Apply Squeeze theorem and (\ref{ln/ln}) to obtain
$$
\lim\limits_{n\rightarrow\infty} \frac{\Phi(\rho n+\rho+1)}{\Phi(\rho n+1)}=1.
$$
Part 1 of the Lemma \ref{lem2} is proved.

The proof of part 2 of the lemma obviously follows from (\ref{estgam}):
$$
\Phi(\rho n+1)\ge \ln{(\rho n+1)}-\frac{1}{\rho n+1}> -\frac{1}{\rho +1}, \quad n\ge 1.
$$
Lemma \ref{lem2} is completely proved.
\end{proof}

\begin{lemma}\label{lemma1}
Let $\rho_0\in(0,1)$ and $t_1>0$. Then series (\ref{derML}) is uniformly convergent with respect to \(t \in [t_1, T]\) and \(\rho \in [\rho_0, 1]\).
\end{lemma}

\begin{proof}
Since $t\in[t_1,T]$ and $\rho\in[\rho_0,1]$ we get
$$
\sum\limits_{n=1}^\infty \left| (-1)^nnt^{\rho n}\frac{\ln{t}-\Phi(\rho n+1)}{\Gamma(\rho n+1)}\right|\le \sum\limits_{n=1}^\infty nT^{n}\frac{|\ln{t}-\Phi(\rho n+1)|}{\Gamma(\rho_0 n+1)}.
$$
Using inequality (\ref{estgam}), we have
$$
\sum\limits_{n=1}^\infty \frac{nT^{n}}{\Gamma(\rho_0 n+1)}\left(\max\limits_{t\in\{t_1,T\}}|\ln{t}|+\ln{(n+1)}+\frac{1}{2}\right)=\sum\limits_{n=1}^\infty z_n,
$$
where
$$
z_n=\frac{nT^{n}}{\Gamma(\rho_0 n+1)}\left(\max\limits_{t\in\{t_1,T\}}|\ln{t}|+\ln{(n+1)}+\frac{1}{2}\right).
$$
Let us consider the following proportion
$$
\frac{z_{n+1}}{z_n}=\frac{(n+1)T^{n+1}\Gamma(\rho_0 n+1)}{nT^{n}\Gamma(\rho_0 n+\rho_0+1)}\frac{\max\limits_{t\in\{t_1,T\}}|\ln{t}|+\ln{(n+2)}+\frac{1}{2}}{\max\limits_{t\in\{t_1,T\}}|\ln{t}|+\ln{(n+1)}+\frac{1}{2}}.
$$
By the equality $\Gamma(\rho n+1)=\rho n\Gamma(\rho n)$, we have
$$
\frac{z_{n+1}}{z_n}=T\frac{\Gamma(\rho_0 n)}{\Gamma(\rho_0 n+\rho_0)}\frac{\max\limits_{t\in\{t_1,T\}}|\ln{t}|+\ln{(n+2)}+\frac{1}{2}}{\max\limits_{t\in\{t_1,T\}}|\ln{t}|+\ln{(n+1)}+\frac{1}{2}}.
$$
Apply Lemma \ref{lem1}, to get
$$
\lim\limits_{n\rightarrow\infty} \frac{z_{n+1}}{z_n}=0.
$$
Thus, by D'Alembert's Ratio Test, we conclude that the series \(\sum\limits_{n=1}^\infty z_n\) is convergent. Therefore, by the Weierstrass M-Test, we deduce that the series (\ref{derML}) is uniformly convergent.

Lemma \ref{lemma1} is proved.
\end{proof}

\begin{lemma}\label{n>n_0}
   Let $0<\rho<1$. Then for $t\in\left(0, \min{\left(\frac{1}{2^{\frac{1}{\rho}}}, \frac{1}{e^{\frac{7}{2}}}\right)}\right]$ the inequalities $y_{n+1}>y_{n},\,\,n\ge 1$, hold.
\end{lemma}
\begin{proof}
Let us consider the following proportion
$$
\frac{y_{n+1}}{y_n}=t^\rho \frac{(n+1)\Gamma(\rho n+1)}{n\Gamma(\rho n+\rho+1)} \frac{\ln(t)-\Phi(\rho n+\rho+1)}{\ln(t)-\Phi(\rho n+1)}.
$$
According to equality $\Gamma(\rho n+1)=\rho n\Gamma(\rho n)$ we have
$$
\frac{y_{n+1}}{y_n}=t^\rho \frac{\Gamma(\rho n)}{\Gamma(\rho n+\rho)} \frac{\ln(t)-\Phi(\rho n+\rho+1)}{\ln(t)-\Phi(\rho n+1)}
$$
\begin{equation}\label{nisbat}
=t^\rho \frac{\Gamma(\rho n)}{\Gamma(\rho n+\rho)} \left(1+\frac{\Phi(\rho n+1)-\Phi(\rho n+\rho+1)}{\ln(t)-\Phi(\rho n+1)}\right)
\end{equation}
Now we will show that
$$
\left|\frac{\Phi(\rho n+1)-\Phi(\rho n+\rho+1)}{\ln(t)-\Phi(\rho n+1)}\right|<1.
$$

Based on estimate (\ref{estgam}), we examine four cases to prove the estimate above.

We can rewrite the expression $\Phi(z)$ according to the estimate (\ref{estpsi}) in the following form:
$$
\Phi(z)=\ln(z)-\frac{\theta(z)}{z}, \quad z\in(0, \infty),
$$
where function $\theta(z)$ is continuous and satisfies $\frac{1}{2}\le\theta(z)\le1$ for all $z\in(0, \infty)$.
From this we get:
$$
\left|\frac{\Phi(\rho n+1)-\Phi(\rho n+\rho+1)}{\ln{t}-\Phi(\rho n+1)}\right|
$$
$$
= \left|\frac{\ln{(\rho n+1)}-\frac{\theta(\rho n+1)}{\rho n+1}-\ln{(\rho n+\rho+1)}+\frac{\theta(\rho n+\rho+1)}{\rho n+\rho+1}}{\ln{t}-\ln{(\rho n+1)}+\frac{\theta(\rho n+1)}{\rho n+1}}\right|
$$
$$
= \left|\frac{\ln{\left(\frac{\rho n+1}{\rho n+\rho+1}\right)}-
\frac{(\rho n+1)(\theta(\rho n+1)-\theta(\rho n+\rho+1))+\rho\theta(\rho n+1)}{(\rho n+1)(\rho n+\rho+1)}}{\ln{t}-\ln{(\rho n+1)}+\frac{\theta(\rho n+1)}{\rho n+1}}\right|
$$
$$
\le \left|\frac{\ln{\left(\frac{\rho n+1}{\rho n+\rho+1}\right)}}{\ln{t}-\ln{(\rho n+1)}+\frac{\theta(\rho n+1)}{\rho n+1}}\right|
+\left|\frac{\frac{(\rho n+1)(\theta(\rho n+1)-\theta(\rho n+\rho+1))+\rho\theta(\rho n+1)}{(\rho n+1)(\rho n+\rho+1)}}{\ln{t}-\ln{(\rho n+1)}+\frac{\theta(\rho n+1)}{\rho n+1}}\right|
$$
$$
\le \frac{1}{|\ln{t}-\ln{(\rho n+1)}|-\frac{\theta(\rho n+1)}{\rho n+1}}+\frac{\left|\frac{\theta(\rho n+1)-\theta(\rho n+\rho+1)}{(\rho n+\rho+1)}+\frac{\rho\theta(\rho n+1)}{(\rho n+1)(\rho n+\rho+1)}\right|}{|\ln{t}-\ln{(\rho n+1)}|-\frac{\theta(\rho n+1)}{\rho n+1}}.
$$
According to $\theta(\rho n+1), \theta(\rho n+\rho+1)\in [1/2, 1]$ we have
$$
\left|\frac{\Phi(\rho n+1)-\Phi(\rho n+\rho+1)}{\ln{t}-\Phi(\rho n+1)}\right|\le\frac{1}{|\ln{t}|-1}+\frac{\frac{1}{2}+1}{|\ln{t}|-1}
=\frac{\frac{5}{2}}{|\ln{t}|-1}.
$$

If $\ln{t}< -\frac{7}{2}$ then we have
\begin{equation}\label{qh}
\left|\frac{\Phi(\rho n+1)-\Phi(\rho n+\rho+1)}{\ln{t}-\Phi(\rho n+1)}\right|<  1.
\end{equation}

Since $\Gamma(\rho n) \le \Gamma(\rho n+\rho)$, the equality (\ref{nisbat}) and the inequality (\ref{qh}) we have the following inequality
$$
\frac{y_{n+1}}{y_n}\le t^\rho  \left(1+\frac{\Phi(\rho n+1)-\Phi(\rho n+\rho+1)}{\ln(t)-\Phi(\rho n+1)}\right)<2t^\rho.
$$
If we select $t\in\left(0, \min{\left(\frac{1}{2^{\frac{1}{\rho}}}, \frac{1}{e^{\frac{7}{2}}}\right)}\right]$, then we have $\frac{y_{n+1}}{y_n}<1$. For such $t$, according to assertion 2 of Lemma \ref{lem2}, the sequence $y_n$, defined by the equality (\ref{sery_n}), is negative for all $n\ge 1$. Therefore, we get the inequality $y_{n+1}>y_{n}$.
Lemma \ref{n>n_0} is proved.
\end{proof}

\begin{theorem}\label{main_lemma1}
Let \(\rho_0 \in (0,1)\). Then, for any \(t \in \left(0, \min\left(\frac{1}{2^{\frac{1}{\rho_0}}}, \frac{1}{e^{\frac{7}{2}}}\right)\right]\), the Mittag-Leffler function \(E_\rho(-t^\rho)\) is monotonically increasing in \(\rho \in [\rho_0, 1]\).
\end{theorem}

\begin{proof} According to Lemma \ref{lemma1} the series (\ref{derML}) converges absolutely. Let us divide this series into groups as follows
$$
\frac{d}{d\rho}E_{\rho}(-t^\rho)=\sum\limits_{n=1}^\infty (-1)^ny_n=-(y_1-y_2)-(y_3-y_4)-\dots
$$
Then by Lemma \ref{n>n_0} it follows that $\frac{d}{d\rho}E_{\rho}(-t^\rho)>0$ for all $t$ and $\rho$ from the conditions of the lemma.

Theorem \ref{main_lemma1} is proved.
\end{proof}

\begin{theorem}\label{main_lemma2}
Let \(\rho_0 \in (0,1)\). Then, for any \(t \in \left(0, \min\left(\frac{1}{2^{\frac{1}{\rho_0}}}, \frac{1}{e^{\frac{13}{6}}}\right)\right]\), the function \(t^{\rho-1}E_{\rho,\rho}(-t^\rho)\) is monotonically decreasing in \(\rho \in [\rho_0, 1]\).
\end{theorem}

The theorem is proved in the same way as Theorem \ref{main_lemma1}.

Using Theorems \ref{main_lemma1} and \ref{main_lemma2}, we can update the results presented in several previous works. In particular, in \cite{AA1, AA2}, the inverse problem of determining the order of the fractional derivative is solved under the over-determination condition
\begin{equation}\label{o}
||u(x, t_0)||_{L_2(\Omega)}^2 = d_0,
\end{equation}
where $t_0$ is a sufficiently large number. Our results show that under the over-determination condition (\ref{o}), the time instant $t_0$ can also be chosen sufficiently small, which is quite useful for real-life applications.

\section{Initial-boundary value problem}

Let us first give the definition of the solution to the problem (\ref{prob1}).
 \begin{definition}\label{def1} A function $u(x,t)$ with the properties
\begin{enumerate}
  \item$u(x,t)\in C(\overline{\Omega}\times [0.T])$,
  \item
  $D_t^\rho u(x,t), \Delta u(x,t)\in C(\overline{\Omega} \times (0.T])$,
\end{enumerate}
and satisfying conditions
(\ref{prob1})  is called \textbf{the solution} of the problem (\ref{prob1}).
\end{definition}

 Now we present theorem about the solution of problem (\ref{prob1}).

\begin{theorem}\label{main1}
 Let function $\varphi (x)\in I_{M}(\Omega) $. Then problem (\ref{prob1}) has a unique solution:
\begin{equation}\label{prob.s}
u(x,t)=\sum\limits_{k=1}^{\infty}\varphi_k E_\rho(-\lambda_k t^\rho )v_k(x),
    \end{equation}
    where $\varphi_k$ are the Fourier coefficients of function $\varphi(x)$.
\end{theorem}

\begin{proof} According to the Fourier method, we will seek the solution to problem (\ref{prob1}) in the form
$$
	u(x,t)=\sum\limits_{k=1}^\infty T_k(t)v_k(x),
$$
where $T_k(t)=(u(\cdot, t),v_k)$ are the Fourier coefficients of the function $u(x, t)$ and are unknown.

Multiply the orthonormal eigenfunctions $v_k(x)$ to the equation in problem (\ref{prob1}), to get
$$
	(D_t^\rho u(\cdot, t),v_k)+(-\Delta u(\cdot, t),v_k)=0.
$$
We have $(D_t^\rho u(\cdot, t),v_k)=D_t^\rho T_k(t)$ and since operator $A$  is self-adjoint, then  $(-\Delta u(\cdot, t),v_k)=\lambda_k T_k(t)$. Therefore, to determine $T_k(t)$ we obtain the following Cauchy problem
$$
D_t^\rho T_k(t)+\lambda_k T_k(t)=0,\quad \quad T_k(0)=\varphi_k.
$$
This problem has a unique solution (see, for example, \cite{Gor}, p. 174):
$$
	T_k(t)=\varphi_k E_{\rho} (-\lambda_k t^\rho).
$$

Let us show that the operators $ (-\Delta)$ and $D_t^\rho $ can be applied term-by-term to series (\ref{prob.s}) and the resulting series converges uniformly in $(x, t) \in  (\overline{\Omega}\times (0, T])$:
   \[
-\Delta u(x,t)=\sum\limits_{k=1}^{\infty} \lambda_k \varphi_k E_\rho(-\lambda_k t^\rho )v_k(x).
\]
   Using Lemma \ref{mll4} and applying the Cauchy-Bunyakovsky inequality we get
\[
|-\Delta u(x,t)|\le\sum\limits_{k=1}^{\infty} |\lambda_k \varphi_k E_\rho(-\lambda_k t^\rho )v_k(x)|\le \sum\limits_{k=1}^{\infty} \lambda_k |\varphi_k| |v_k(x)| \frac{C}{1+\lambda_kt^\rho}
\]
$$
\le Ct^{-\rho}\sum\limits_{k=1}^{\infty} |\varphi_k| |v_k(x)|=Ct^{-\rho}\sum\limits_{k=1}^{\infty} |\varphi_k|(\sqrt{\lambda_k})^{\left(\left[\frac{N}{2}\right]+1\right)}(\sqrt{\lambda_k})^{-\left(\left[\frac{N}{2}\right]+1\right)}|v_k(x)|
$$
$$
\le Ct^{-\rho}\left(\sum\limits_{k=1}^{\infty} |\varphi_k|^2\lambda_k^{\left(\left[\frac{N}{2}\right]+1\right)}\right)^{\frac{1}{2}}\left(\sum\limits_{k=1}^{\infty}\lambda_k^{-\left(\left[\frac{N}{2}\right]+1\right)}v_k^2(x)\right)^{\frac{1}{2}}
$$

Therefore, if $\varphi(x)\in I_M(\Omega)$, then the series $\sum\limits_{k=1}^{\infty} |\varphi_k|^2\lambda_k^{\left(\left[\frac{N}{2}\right]+1\right)}$ converges.
Furthermore, according to Lemma \ref{Illem1} the series $\sum\limits_{k=1}^{\infty}\lambda_k^{-\left(\left[\frac{N}{2}\right]+1\right)}v_k^2(x)$ converges uniformly in a closed domain $\overline{\Omega}$. Hence, $ -\Delta u(x, t) \in C(\overline{\Omega}\times (0, T])$.

   From equation (\ref{prob1}) one has
$D_t^\rho u(x, t) = \Delta u(x, t)$, $t > 0$, and hence we get $D_t^\rho u(x, t) \in C(\overline{\Omega}\times (0, T])$.

The uniqueness of the solution is proved in the standard way (see, for example, \cite{AU}).
\end{proof}

\section{Inverse problem}

Before moving on to the inverse problem, let us examine the uniform convergence of the following series
\begin{equation}\label{ser}
    \sum\limits_{k=1}^\infty \varphi_kv_k(x) \frac{d}{d\rho}E_{\rho}(-\lambda_k t^\rho),
\end{equation}
where $\varphi_k=(\varphi,v_k)$. One has

\begin{lemma}\label{Gruncon}
    Let $\rho_0\in(0,1)$ and $t_1>0$. If the function $\varphi(x)\in I_M(\Omega)$. Then the series (\ref{ser}) is uniformly convergent with respect to $x\in\overline{\Omega},\,\, t\in[t_1,T],\,\, \rho\in[\rho_0,1]$.
\end{lemma}
\begin{proof} In Lemma \ref{lemma1} the uniform convergence for all $t\in[t_1,T]$ and $\rho\in[\rho_0,1]$ of the series (\ref{derML}) for $\frac{d}{d\rho}E_{\rho}(-\lambda_k t^\rho)$ is proved. Therefore there exists a number $C_0$ such that $\left|\frac{d}{d\rho}E_{\rho}(-\lambda_k t^\rho)\right|\le C_0$. Next, the uniform convergence of the series (\ref{ser})  is proved using the same arguments as in the proof of Theorem \ref{main1}.
Lemma \ref{Gruncon} is proved.
\end{proof}

From Lemma \ref{Gruncon} we get the following corollary.

\begin{corollary}\label{corollary1}
    Let $\rho_0\in(0,1)$ and $t_1>0$. Then for any number $\varepsilon>0$, there exists $n_0=n_0(\varepsilon, M, \rho_0, t_1)$ such that the following inequality holds for all $x\in\overline{\Omega},\,\, t\in [t_1,T], \rho\in[\rho_0,1]$ and $\varphi(x)\in I_M(\Omega)$:
    $$
    \sum\limits_{k=n_0+1}^\infty  \left|\varphi_kv_k(x) \frac{d}{d\rho}E_{\rho}(-\lambda_k t^\rho)\right|<\varepsilon.
    $$
\end{corollary}

Now we move on to the study of the inverse problem. Let $x_0\in\Omega$ be fixed and we have the measured data at one time instant $t_0>0$ given as
\begin{equation}\label{adcon}
    u(x_0,t_0)=d_0,
\end{equation}
where $d_0$ is a fixed number. If we introduce the following nonlinear function
\begin{equation}\label{adfunc}
    G(\rho)=\sum\limits_{k=1}^\infty \varphi_k E_{\rho}(-\lambda_k t_0^\rho)v_k(x_0),
\end{equation}
 then the equation (\ref{adcon}) for determining the unknown parameter $\rho$ can be written as $G(\rho)=d_0$.

 Based on Corollary \ref{corollary1}, we introduce the following definition of the solution of the Inverse Problem (\ref{prob1}), (\ref{adcon}).
  \begin{definition}\label{def2} The pair $\{u(x,t), \rho\}$ of the function $u(x,t)$ and the parameter $\rho$ with the properties: $\rho\in [\rho_0, 1]$, $\rho_0\in (0, 1)$, $u(x,t)$ satisfies the conditions of Definition \ref{def1} and together with $\rho$ satisfies the over-determination condition (\ref{adcon}) is called a \textbf{solution} of the Inverse Problem (\ref{prob1}), (\ref{adcon}).
\end{definition}

Now we present the theorem about the solution of Inverse Problem (\ref{prob1}), (\ref{adcon}).

\begin{theorem}\label{Main_teorem}
Let $\varepsilon > 0$ be a fixed number and
\( t_0 \in \left(0, \min\left(\frac{1}{2^{1/\rho_0}}, \frac{1}{e^{7/2}}\right)\right] \).
Denote \( n_0 = n_0(\varepsilon, M, \rho_0, t_0) \) as in Corollary \ref{corollary1}. Suppose that for some $x_0 \in \Omega$ the quantities $v_k(x)$ and $\varphi_k$ of $\varphi(x)\in I_M(\Omega)$ satisfy the following conditions: \begin{enumerate}
    \item $\varphi_k \ge 0$ and $v_k(x_0) \ge 0$ for $k = 1, 2, \ldots, n_0$;
    \item $\varphi_{k_0}$ and $v_{k_0}(x_0)$ satisfying following inequality for some $k_0\in\{1, 2, \ldots, n_0\}$
    $$
    \varphi_{k_0}v_{k_0}(x_0)>\frac{\varepsilon}{M_{k_0}},
    $$
    where $M_{k_0}=\frac{d}{d\rho}E_{\rho}(-\lambda_{k_0} t_0^\rho)$.
\end{enumerate}
Then the nonlinear function $G(\rho)$, defined by (\ref{adfunc}), is strictly monotonically increasing on $\rho \in [\rho_0,1]$. Moreover, Inverse Problem (\ref{prob1}), (\ref{adcon}) has a unique solution if and only if
\begin{equation}\label{d0}
   d_0\in [G(\rho_0), G(1)].
\end{equation}
\end{theorem}

\begin{proof} First, we will show that the function $G(\rho)$ is strictly monotonic. To demonstrate the monotonicity of the function $ G(\rho) $, it is sufficient to show that for any fixed $\rho_1, \rho_2 \in [\rho_0, 1]$ with $\rho_2 > \rho_1$, the inequality $ G(\rho_2) - G(\rho_1) > 0 $ holds.
So
$$
G(\rho_2)-G(\rho_1)=\sum\limits_{k=1}^\infty \varphi_k v_k(x_0) [E_{\rho_2}(-\lambda_kt_0^{\rho_2})-E_{\rho_1}(-\lambda_k t_0^{\rho_1})].
$$
Let $\varepsilon>0$ be an arbitrary given number, and $n_0$ be a number from Corollary \ref{corollary1}. Then we may write
$$
G(\rho_2)-G(\rho_1)=\sum\limits_{k=1}^{n_0} \varphi_k v_k(x_0) [E_{\rho_2}(-\lambda_kt_0^{\rho_2})-E_{\rho_1}(-\lambda_k t_0^{\rho_1})]
$$
$$
+\sum\limits_{k=n_0+1}^{\infty} \varphi_k v_k(x_0) [E_{\rho_2}(-\lambda_kt_0^{\rho_2})-E_{\rho_1}(-\lambda_k t_0^{\rho_1})]=I_1+I_2.
$$
Now we estimate the following expression
$$
|I_2|=\left|\sum\limits_{k=n_0+1}^{\infty} \varphi_k v_k(x_0) [E_{\rho_2}(-\lambda_kt_0^{\rho_2})-E_{\rho_1}(-\lambda_k t_0^{\rho_1})]\right|
$$
$$
\le \sum\limits_{k=n_0+1}^{\infty}\left|\varphi_k v_k(x_0) [E_{\rho_2}(-\lambda_kt_0^{\rho_2})-E_{\rho_1}(-\lambda_k t_0^{\rho_1})]\right|.
$$
Since $ E_\rho(-\lambda_k t_0^\rho) $ is continuously differentiable with respect to $\rho$ at each point $\rho>0$, by Lagrange theorem, there exists $q_1=q_1(k)\in [\rho_1, \rho_2]$ such that ($k>n_0$)
$$
E_{\rho_2}(-\lambda_k t_0^{\rho_2}) - E_{\rho_1}(-\lambda_k t_0^{\rho_1}) = \frac{d}{d\rho}(E_{\rho}(-\lambda_k t_0^{\rho}))\bigg|_{\rho=q_1} (\rho_2 - \rho_1)
$$
$$
:=\frac{d}{d\rho}(E_{q_1}(-\lambda_k t_0^{q_1})) (\rho_2 - \rho_1).
$$
From this we get:
$$
|I_2|\le \sum\limits_{k=n_0+1}^{\infty}\left|\varphi_k v_k(x_0)\frac{d}{d\rho}E_{q_1}(-\lambda_k t_0^{q_1})\right| (\rho_2 - \rho_1)
$$
$$
=(\rho_2-\rho_1)\sum\limits_{k=n_0+1}^{\infty}\left|\varphi_k v_k(x_0)\frac{d}{d\rho}E_{q_1}(-\lambda_k t_0^{q_1})\right|.
$$
Apply Corollary \ref{corollary1} to obtain:
$$
|I_2|\le (\rho_2-\rho_1)\varepsilon.
$$

Let us consider $I_1$:
$$
I_1=\sum\limits_{k=1}^{n_0} \varphi_k v_k(x_0) [E_{\rho_2}(-\lambda_kt_0^{\rho_2})-E_{\rho_1}(-\lambda_k t_0^{\rho_1})].
$$
Again, by Lagrange theorem, there exists $q_2=q_2(k)\in [\rho_1, \rho_2]$ such that ($k\leq n_0$)
$$
E_{\rho_2}(-\lambda_k t_0^{\rho_2}) - E_{\rho_1}(-\lambda_k t_0^{\rho_1}) = \frac{d}{d\rho}(E_{\rho}(-\lambda_k t_0^{\rho}))\bigg|_{\rho=q_2} (\rho_2 - \rho_1)
$$
$$
:=\frac{d}{d\rho}(E_{q_2}(-\lambda_k t_0^{q_2})) (\rho_2 - \rho_1).
$$
From this and the fact that $ \frac{d}{d\rho}E_{q_2}(-\lambda_k t_0^{q_2}) = M_k $, we have
$$
I_1= (\rho_2-\rho_1) \sum_{k=1}^{n_0} \varphi_k v_k(x_0)M_k.
$$
If $\varphi_k \ge 0$ and $v_k(x_0) \ge 0$ for $k = 1, 2, \ldots, n_0$, and for some $k_0\in (1, 2, \ldots, n_0)$, $\varphi_{k_0} > 0$ and $v_{k_0}(x_0) > 0$, then we get
$$
I_1>(\rho_2-\rho_1)\varphi_{k_0}v_{k_0}(x_0)M_{k_0}.
$$
Since $\varphi_{k_0}v_{k_0}(x_0) > \frac{\varepsilon}{M_{k_0}}$ we have $G(\rho_2) - G(\rho_1) > 0$. Therefore, $G(\rho)$ is monotonically increasing.

Obviously, if condition (\ref{d0}) is not met, then the solution to the Inverse Problem does not exist.

Theorem is proved.
\end{proof}

\noindent\textbf{Remark 1.}\label{rem}
{\it Let $M_0 = \min\limits_{1 \leq k \leq n_0, \rho\in [\rho_0, 1]} M_k(\rho)$. Then from the given proof it follows that condition (2) of Theorem \ref{Main_teorem} can be replaced by the condition}
\begin{equation}\label{new}
   \sum_{k=1}^{n_0} \varphi_k v_k(x_0)> \frac{\varepsilon}{M_{0}}.
\end{equation}

Note that in Theorem 2 of the paper \cite{GongLi}  by G. Li, Z. Wang, X. Jia, Y. Zhang, in particular, condition (2) is missing. From the proof of this theorem it becomes obvious that if condition (2) (or the condition (\ref{new}) in Remark \ref{rem}) is missing, then it is easy to construct a counterexample showing that the assertion of the theorem is not satisfied.

\section{An example of Sh.A. Alimov}
It should be noted that the formulation of the inverse problem considered in the work \cite{GongLi}  by G. Li, Z. Wang, X. Jia, Y. Zhang  is interesting, but Theorem \ref{Main_teorem} requires rather strict conditions on the corresponding eigenfunctions and coefficients of the initial function. In this connection, a natural question arises: does an initial-boundary value problem for the subdiffusion equation exist for which the inverse problem with condition (\ref{adcon}) has a unique solution? The following simple example by Sh.A. Alimov provides a positive answer to this question.

Let $A$ be the operator $Av(x,y)=-\Delta v(x,y)$ with boundary conditions
\begin{equation}\label{1}
\frac{\partial v}{\partial \nu}  = 0, \quad (x,y)\in \partial\Omega\times (0,H),
\end{equation}
and
\begin{equation}\label{2}
\frac{\partial v}{\partial y}(x,0) = 0, \quad \frac{\partial v}{\partial y}(x,H) + hv(x,H) = 0, \quad x\in \Omega.
\end{equation}
Here $\Delta =\sum\limits_{k=1}^N \frac{\partial^2}{\partial x_k^2}+\frac{\partial^2}{\partial y^2}$ is the Laplace operator, $\nu$ is an outer normal to $\Omega$ and $h>0$.

Condition (\ref{1}) means that the boundary of domain $\Omega$ is reliably isolated. Condition (\ref{2}) means that penetration and diffusion is possible only through one of the bases of the cylindrical region $D=\Omega\times (0,H)$, and the greater the concentration on this base, the greater the penetration flow.

Consider the Cauchy type problem:
\begin{equation}\label{3}
\left\{
\begin{aligned}
&\partial_t^{\rho} u(x,y,t) + Au(x,y,t) = 0, \quad 0<t\leq T,\quad (x,y)\in D, \\
& \lim\limits_{t\rightarrow 0} \partial^{\rho-1}u(x,y,t) = \phi(x,y),\quad (x,y)\in D,
\end{aligned}
\right.
\end{equation}
where $\phi$ is a given continuous function and $\partial^\rho_x$ is the fractional Riemann-Liouville derivative defined as (see, for example, \cite{Klb}, p. 70):
$$
\partial_x^\rho u(x)=\frac{1}{\Gamma
(1-\rho)}\frac{d}{dx}\int\limits_0^x\frac{u(\xi)}{(x-\xi)^{\rho}} d\xi,
\quad x>0.
$$

The first eigenvalue of operator $A$ is negative: $\lambda_1=- \mu^2$, where $\mu = \mu(hH)$ is the unique positive solution to the equation
$$
\mu\tanh \mu\ =\ hH,
$$
and the corresponding first eigenfunction is positive: $u_1(x,y)>0$.

Further reasoning is valid for an arbitrary initial-boundary value problem (\ref{3}), in which the operator $A$ has these two properties, i.e.
$$
\lambda_1 < 0, \quad \lambda_1<\lambda_2\leq  \lambda_3\leq\cdots
$$
and
\begin{equation}\label{eigenfunction}
u_1(x,y)>0.
\end{equation}

Under certain conditions on the initial function, the unique solution to problem (\ref{3}) has the form (see Theorem \ref{main1}):
\begin{equation}\label{sol}
u(x,y,t)\ = \sum\limits_{m=1}^\infty t^{\rho-1} E_{\rho,\rho}(-\lambda_m t^\rho)(\phi,u_m)u_m(x,y),
\end{equation}
which absolutely and uniformly converges on $(x, y)\in \overline{D}$
 for each $t\in (0, T]$.

 Let us consider the inverse problem with the over-determination condition:
\begin{equation}\label{ad}
U(\rho; t_0)\equiv u(x_0, y_0, t_0) =d_0,
\end{equation}
at a fixed time instant $t_0$ and a point $(x_0, y_0)\in D$.

One has

\begin{lemma}\label{U}
Let $(x_0, y_0)$ be a completely arbitrary point in domain $D$ and the Fourier coefficient $\phi_1=(\phi,u_1)$ of $\phi(x,y)$ is not zero. Then there is a number $T_0=T_0(\rho, \lambda^*)$, $\lambda^* = \min\limits_{j} |\lambda_j|$, such that for all $t_0\geq T_0$ the function $U(\rho; t_0)$ is strictly  monotone: if $\phi_1 < 0$, then $U(\rho; t_0)$ increases, and if
$\phi_1 > 0$, then  $U(\rho; t_0)$ decreases as a function of $\rho$.
\end{lemma}

Naturally, equation (\ref{ad}) with respect to $\rho$ does not have a solution for every $d_0$: a necessary and sufficient condition for the existence of a solution to this equation is the fulfillment of inequality
\begin{equation}\label{cond0}
    \inf\limits_{\rho\in (0, 1)}U(\rho, t_0)<d_0<\max\limits_{\rho\in (0, 1)}U(\rho, t_0)
\end{equation}

The following theorem gives a positive answer to the above  question.

\begin{theorem}\label{ip1}Suppose the conditions of Lemma \ref{U} are satisfied. Let $d_0$ satisfy condition (\ref{cond0}) and let $t_0\geq T_0$. Then inverse problem (\ref{3}), (\ref{ad})  has a unique solution $\{u(x,y,t), \rho \}$.
\end{theorem}

First, we present a proof of Lemma \ref{U}.
\begin{proof}
By virtue of formula (\ref{sol}), function $U(\rho; t_0) $ has the form
\begin{equation}\label{Ut0}
U(\rho; t_0)\ = \sum\limits_{m=1}^\infty t_0^{\rho-1} E_{\rho,\rho}(-\lambda_m t_0^\rho)\phi_m u_m(x_0,y_0).
\end{equation}

Let us recall the asymptotic formulas (see \cite{Dzh66}, p. 134)
\begin{equation}\label{>0}
    E_{\rho, \rho} (z)=\frac{1}{\rho} z^{\frac{1}{\rho} -1}e^{z^{\frac{1}{\rho}}}+ O\bigg(\frac{1}{z}\bigg), \quad z>1,
\end{equation}
\begin{equation}\label{<0}
 E_{\rho, \rho} (-r)=\frac{r^{-2}}{\Gamma(-\rho)} + O(r^{-3}), \quad r>1.
\end{equation}
Taking into account the properties of eigenvalues, we conclude that in the sum (\ref{Ut0}) in the first few terms the argument of the function $E_{\rho,\rho}$ is positive, and then all the rest are negative. Therefore, using estimates (\ref{>0}) and (\ref{<0}), it is not hard  to see that for $t^\rho_0\lambda^*>1$, $\lambda^* = \min\limits_{j} |\lambda_j|$,
$$
    U(\rho; t_0)= \phi_1 u_1(x_0, y_0) \ \frac{1}{\rho} |\lambda_1|^{\frac{1}{\rho} -1}e^{t_0|\lambda_1|^{\frac{1}{\rho}}}\bigg(1+ O(1) e^{-\epsilon_0 t_0^\rho}\bigg),  \,\,\, \epsilon_0=\lambda_2-\lambda_1>0.
$$
Since function $f(z)=z|\lambda_1|^ze^{t_0|\lambda_1|^z} $ is increasing and eigenfunction $u_1(x, y)$ is positive (see (\ref{eigenfunction})), then function $U(\rho; t_0)$ increases if $\phi_1<0$ and decreases if $\phi_1>0$. Lemma \ref{U} is proved.
\end{proof}

Theorem \ref{ip1} is a consequence of this lemma.

\section{On a result of A.V. Pskhu}

Let us note an interesting result by A.V. Pskhu \cite{Pskhu}, published in 2002, which, as far as we know, represents the first inverse problem for determining the order of fractional derivatives.

Consider the following equation
\begin{equation}\label{Pskhu1}
\partial^\rho_x u(x)-\lambda u(x)=0,
\end{equation}
where, $x\in (0,a],\,\, \rho\in(0,1],\,\, \lambda\in\mathbb{R}$ and $\partial^\rho_x$ is the fractional Riemann-Liouville derivative.

\textbf{Inverse problem for equation (\ref{Pskhu1}).} Find the value of order $\rho,\,\, \rho\in(0, 1]$, and solution $u(x)$, of equation (\ref{Pskhu1}) which obeys conditions
\begin{equation}\label{Pskhu2}
    u(x_0)=u_0,\quad u(x_1)=u_1,
\end{equation}
where $x_0,\,\, x_1\in(0,a],\,\, u_0,\,\, u_1$ are given real numbers.

As the author proved, the problem has a unique solution if $x_0\neq x_1$ and
$$
\frac{x_0}{x_1}< \frac{u_1}{u_0} \leq e^{\lambda (x_1 - x_0)}.
$$

We adopted a slightly different approach to the inverse problem (\ref{Pskhu1}), (\ref{Pskhu2}) and derived the following result.

Let $\rho\in(0,1]$. Consider the following Cauchy problem
\begin{equation}\label{Pskhuprob}
\left\{
\begin{aligned}
&\partial^\rho_x u(x)-\lambda u(x)=0,\quad x\in(0,1],\\
&\lim\limits_{x\rightarrow 0}\partial_x^{\rho-1} u(x)=\varphi,
\end{aligned}
\right.
\end{equation}
where, $\varphi, \lambda\in\mathbb{R}$ are given numbers.

\textbf{Inverse problem for problem (\ref{Pskhuprob}).} Find  the value of order $\rho,\,\, \rho\in(0, 1]$, and solution $u(x)$, of problem (\ref{Pskhuprob}) which obeys condition
\begin{equation}\label{Pskhu3}
    u(x_0)=u_0,
\end{equation}
where $x_0\in(0,1]$ and $u_0$ are given real number.

Note the unique solution to problem (\ref{Pskhuprob})  has the form (see, for example, \cite{Klb}, p. 224)
\begin{equation}\label{solPskhuprob}
    u(x)=\varphi x^{\rho-1} E_{\rho,\rho}(\lambda x^\rho).
\end{equation}

\begin{lemma}\label{Pskhulem}
Let \(x \in (0, 1]\) and $\lambda>0$. The function \(x^{\rho-1}E_{\rho,\rho}(\lambda x^\rho)\) decreases monotonically with respect to \(\rho \in (0, 1]\).
\end{lemma}

\begin{proof} Let $\rho_1, \rho_2\in(0,1],\,\, \rho_1<\rho_2$. We have
$$
x^{\rho_1-1}E_{\rho_1,\rho_1}(\lambda x^{\rho_1})-x^{\rho_2-1}E_{\rho_2,\rho_2}(\lambda x^{\rho_2})=\sum\limits_{k=0}^\infty \frac{\lambda^k x^{\rho_1(k+1)-1}}{\Gamma(\rho_1(k+1))}
$$
$$
-\sum\limits_{k=0}^\infty \frac{\lambda^k x^{\rho_2(k+1)-1}}{\Gamma(\rho_2(k+1))}=\sum\limits_{k=0}^\infty \lambda^kx^{-1}\left(\frac{ x^{\rho_1(k+1)}}{\Gamma(\rho_1(k+1))}-\frac{ x^{\rho_1(k+1)}}{\Gamma(\rho_2(k+1))}\right).
$$
Now consider the following proportion
$$
Q=\frac{\frac{ x^{\rho_2(k+1)}}{\Gamma(\rho_2(k+1))}}{\frac{ x^{\rho_1(k+1)}}{\Gamma(\rho_1(k+1))}}=x^{(\rho_2-\rho_1)(k+1)}\frac{\Gamma(\rho_1(k+1))}{\Gamma(\rho_2(k+1))}.
$$
Since $\Gamma(\rho_1(k+1)) \le \Gamma(\rho_2(k+1))$ and $x\in(0,1]$ we have $Q\le 1$ for all $k\ge 0$. Therefore, $x^{\rho_1-1}E_{\rho_1,\rho_1}(\lambda x^{\rho_1})>x^{\rho_2-1}E_{\rho_2,\rho_2}(\lambda x^{\rho_2})$ for all $\rho_1, \rho_2\in(0,1],\,\, \rho_1<\rho_2$.

Lemma \ref{Pskhulem} is proved.
\end{proof}

The following result is true.
\begin{theorem}\label{Pskhutheorem1}
    Let $x_0\in(0,1]$ and $\lambda>0$. Then the inverse problem (\ref{Pskhuprob}), (\ref{Pskhu3}) has a unique solution $\{u(x), \rho\}$ if and only if
    $$
e^{\lambda x_0} \le \frac{u_0}{\varphi} < +\infty   .
    $$
\end{theorem}
\begin{proof} Using (\ref{solPskhuprob}) and equality (\ref{Pskhu3}) we get
$$
\varphi x_0^{\rho-1}E_{\rho,\rho}(\lambda x_0^\rho)=u_0.
$$
Obviously, according to Lemma \ref{Pskhulem}, the inverse problem (\ref{Pskhuprob}), (\ref{Pskhu3}) has a unique solution.

Theorem \ref{Pskhutheorem1} is proved.
\end{proof}

\begin{theorem}\label{Pskhutheorem}
    Let \(\rho_0 \in (0,1)\), \(x_0 \in \left(0, \min\left(\frac{1}{2^{\frac{1}{\rho_0}}}, \frac{1}{e^{\frac{13}{6}}}\right)\right]\) and $\lambda<0$. Then the inverse problem (\ref{Pskhuprob}), (\ref{Pskhu3}) has a unique solution $\{u(x), \rho\}$ if and only if
    $$
e^{\lambda x_0} \le \frac{u_0}{\varphi} < +\infty   .
    $$
\end{theorem}

The proof of this theorem follows from Theorem \ref{main_lemma2}.

\section*{Acknowledgement}
The authors are grateful to  Sh. A. Alimov for discussions of these results.

The first author acknowledge financial support from the Innovative Development Agency under the Ministry of Higher Education, Science and Innovation of the Republic of Uzbekistan, Grant No F-FA-2021-424.

\bibliographystyle{amsplain}

\end{document}